\pdfoutput=1
\documentclass[11pt, oneside]{article}
\usepackage{geometry} 
\usepackage{graphicx}
\usepackage{caption}
\usepackage{subcaption}
\usepackage{tabularx}
\usepackage{booktabs}

\usepackage{amssymb}
\usepackage{url}
\usepackage{hyperref}
\usepackage{amsmath}

\title{Cabinet of curiosities: the interesting geometry of the angle $\beta = \arccos \left( \left(3 \phi - 1\right)/4 \right)$}
\author{Fang Fang, Klee Irwin\footnote{Corresponding author: \href{mailto:klee@quantumgravityresearch.org}{klee@quantumgravityresearch.org}}, Julio Kovacs, Garrett Sadler\\\emph{Quantum Gravity Research}, Topanga, CA, USA}

\begin{document}

\maketitle

\begin{abstract}
In this paper we present the construction of several aggregates of tetrahedra. Each construction is obtained by performing rotations on an initial set of tetrahedra that either (1) contains gaps between adjacent tetrahedra, or (2) exhibits an aperiodic nature. Following this rotation, gaps of the former case are ``closed'' (in the sense that faces of adjacent tetrahedra are brought into contact to form a ``face junction'') while translational and rotational symmetries are obtained in the latter case. In all cases, an angular displacement of $\beta = \arccos \left(3\phi-1\right)/4$ (or a closely related angle), where $\phi = \left(1+\sqrt{5}\right)/2$ is the golden ratio, is observed between faces of a junction. Additionally, the overall number of \emph{plane classes}, defined as the number of distinct facial orientations in the collection of tetrahedra, is reduced following the transformation. Finally, we present several ``curiosities'' involving the structures discussed here with the goal of inspiring the reader's interest in constructions of this nature and their attending, interesting properties.
\end{abstract}

\section{Introduction}

The present document introduces the reader to the angle $\beta = \arccos \left( \left(3 \phi - 1 \right)/4\right)$, where $\phi = \left(1+\sqrt{5}\right)/2$ is the golden ratio, and its involvement, most notably, in the construction of several interesting aggregates of regular tetrahedra. In the sections below, we will perform geometric rotations on tetrahedra arranged about a common central point, common vertex,  common edge, as well as those of a linear, helical arrangement known as the Boerdijk-Coxeter helix (tetrahelix) \cite{coxeter1974} \cite{boerdijk1952}. In each of these transformations, the angle $\beta$ above appears in the projections of coincident tetrahedral faces. Noteworthy about these transformations is that they have a tendency to bring previously separated faces of ``adjacent'' tetrahedra into contact and to impart a periodic nature to previously aperiodic structures. Additionally, after performing the rotations described below, one observes a reduction in the total number of \emph{plane classes}, defined as the total number of distinct facial  or planar orientations in a given aggregation of polyhedra.

\section{Aggregates of tetrahedra}\label{S:aggregates}

In this section we describe the construction of several interesting aggregates of regular tetrahedra. The aggregates of Sections~\ref{S:comedge} and \ref{S:comvertex} initially contain gaps of various sizes. By performing special rotations of these tetrahedra these gaps are ``closed'' (in the sense that faces of adjacent tetrahedra are made to touch), and, in each case, the resulting angular displacement between coincident faces is either identically equal to $\beta$ or is closely related. In Section~\ref{S:helix}, a rotation by $\beta$ is imparted to tetrahedra arranged in a helical fashion in order to introduce a periodic structure and previously unpossessed symmetries.

\subsection{Aggregates about a common edge}\label{S:comedge}

\begin{figure}[t]
	\centering
	\begin{subfigure}[t]{0.3\textwidth}
		\centering
		\includegraphics[width=\textwidth]{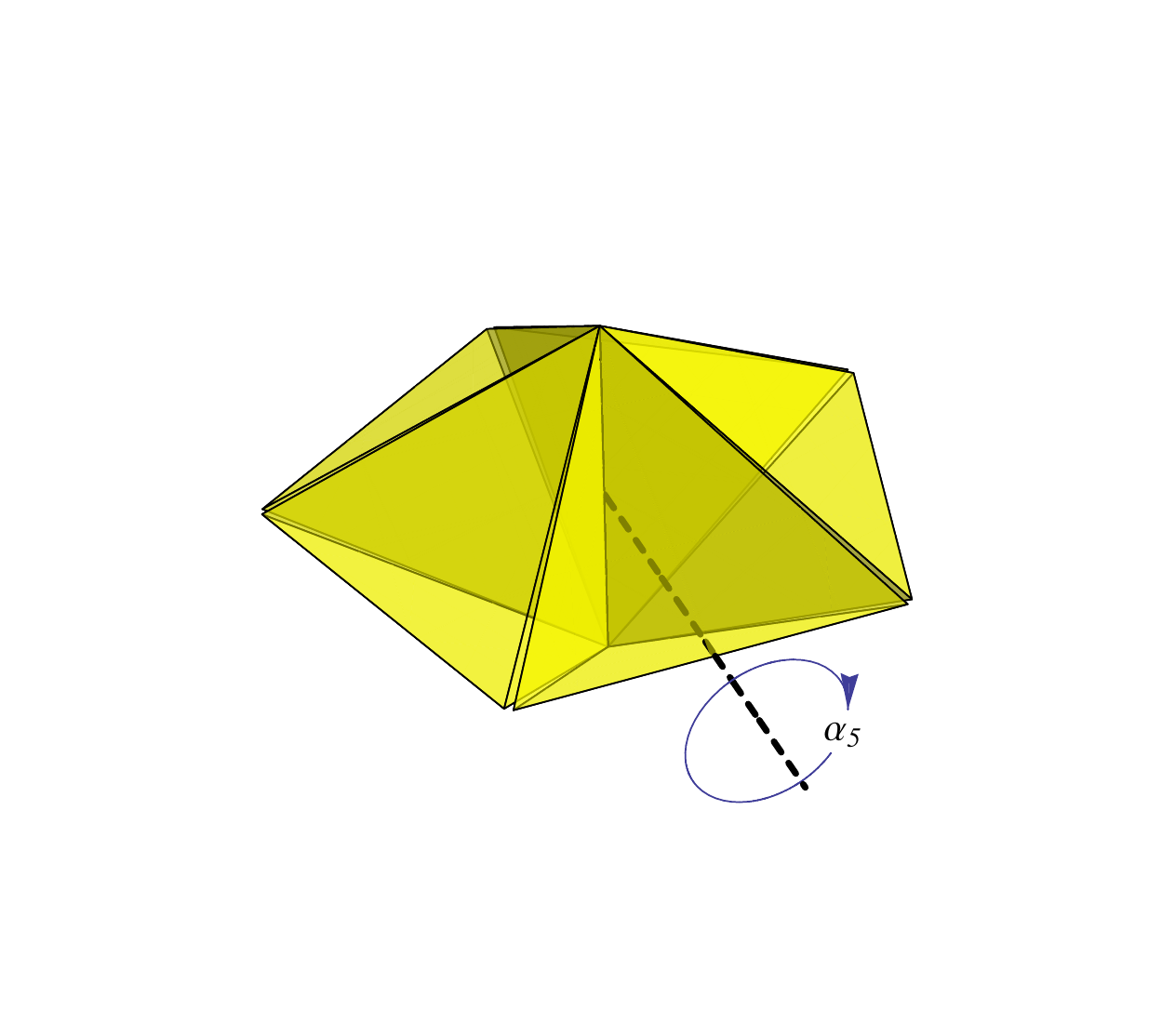}
		\caption{Five tetrahedra arranged about a common edge. In this arrangement, small gaps exist between the faces of adjacent tetrahedra. Each tetrahedron is to be rotated by $\alpha_5$ about an axis passing between the midpoints of its central and peripheral edges.}
		\label{F:fgA}
	\end{subfigure}
	\quad
	\begin{subfigure}[t]{0.3\textwidth}
		\centering
		\includegraphics[width=\textwidth]{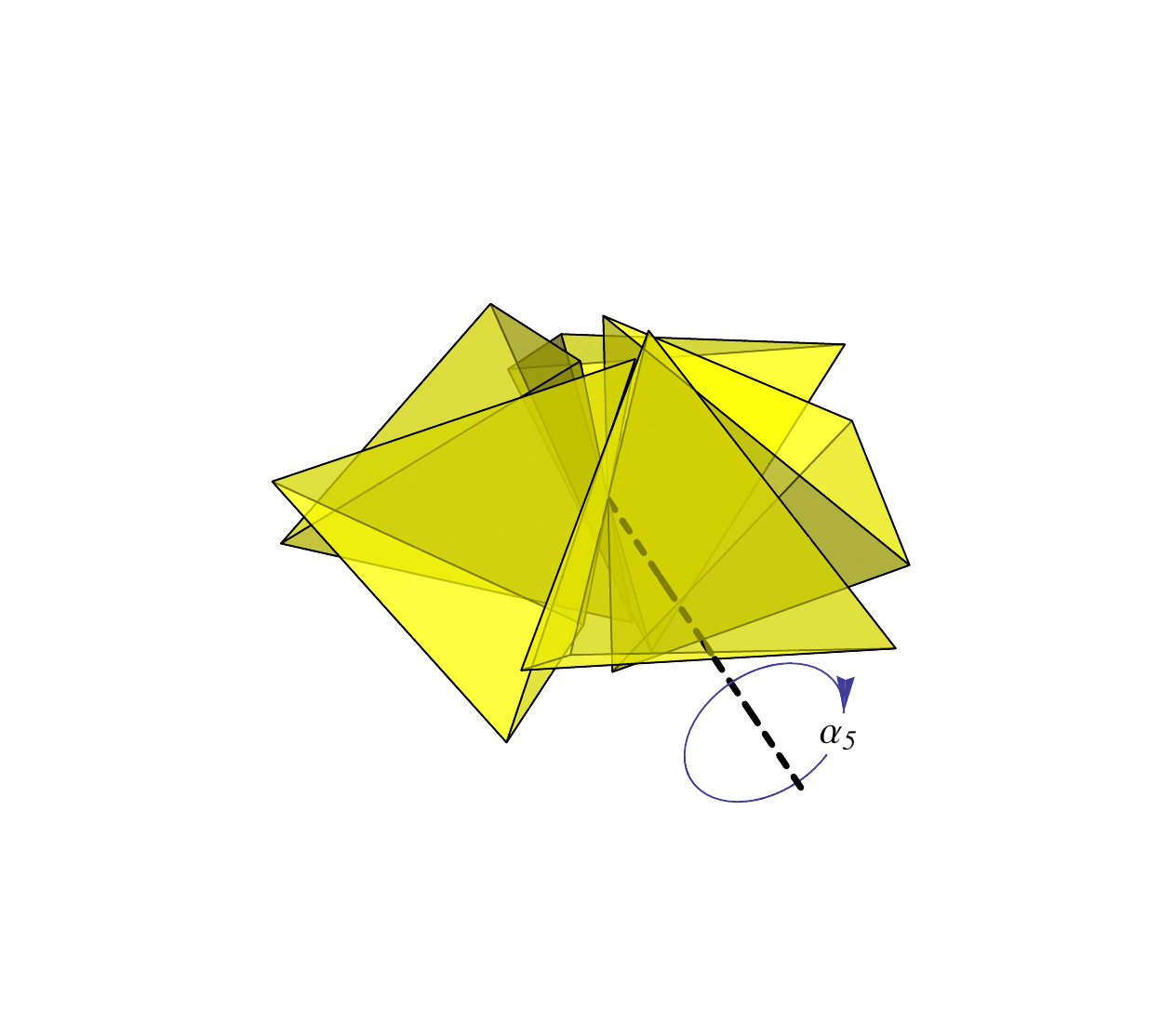}
		\caption{The tetrahedra after rotation. In this arrangement, the faces of adjacent tetrahedra have been brought into contact with one another.}
		\label{F:fgB}
	\end{subfigure}
	\quad
	\begin{subfigure}[t]{0.3\textwidth}
		\centering
		\includegraphics[width=\textwidth]{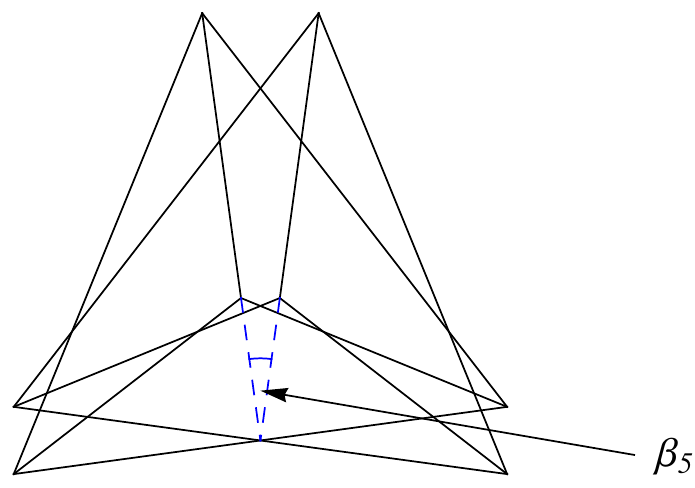}
		\caption{A projection of a ``face junction'' between two coincident faces in \ref{F:fgB}. The angular displacement between the faces is $\beta_5 = \beta$.}
		\label{F:fgC}
	\end{subfigure}
	\caption{``Twisting'' tetrahedra centered about a common central edge to close up gaps between adjacent tetrahedra. When this operation is performed, the angle $\beta = \arccos \left(\left(3\phi - 1\right)/4 \right)$ is produced in the projection of a ``face junction.''}
	\label{F:fivegroup}
\end{figure}

Consider aggregates of \emph{n} regular tetrahedra, $3 \leq n \leq 5$, arranged about a common edge (so that an angle of $2 \pi / n$ is subtended between adjacent tetrahedral centers, see Figure~\ref{F:fgA} for an example with five tetrahedra). In each of these structures, gaps exist between tetrahedra that may be ``closed'' (i.e., faces are made to touch) by performing a rotation of each tetrahedron about an axis passing between the midpoints of its central and peripheral edges through an angle given by
\newline
\begin{equation}\label{E:alphan}
	\alpha_n = \arctan \left( \frac{\sqrt{\cos^2 \frac{\gamma}{2} - \cos^2 \frac{\theta_n}{2}}}{\sin \frac{\gamma}{2} \cos \frac{\theta_n}{2}} \right),
\end{equation}
\newline
where $\gamma = \arccos(1/3)$ is the tetrahedral dihedral angle and $\theta_n = 2\pi/n$. When this is done, an angle, $\beta_n$, is established in the ``face junction'' between coincident pairs of faces such that
\newline
\begin{equation}\label{E:betan}
	\beta_n = 2 \arctan \left( \frac{\sqrt{\cos^2 \frac{\gamma}{2} - \cos^2 \frac{\theta_n}{2}}}{\cos \frac{\theta_n}{2}} \right).
\end{equation}
\newline

\begin{figure}[t]
	\label{F:tg}
	\centering
	\begin{subfigure}[t]{0.3\textwidth}
		\centering
		\includegraphics[width=\textwidth]{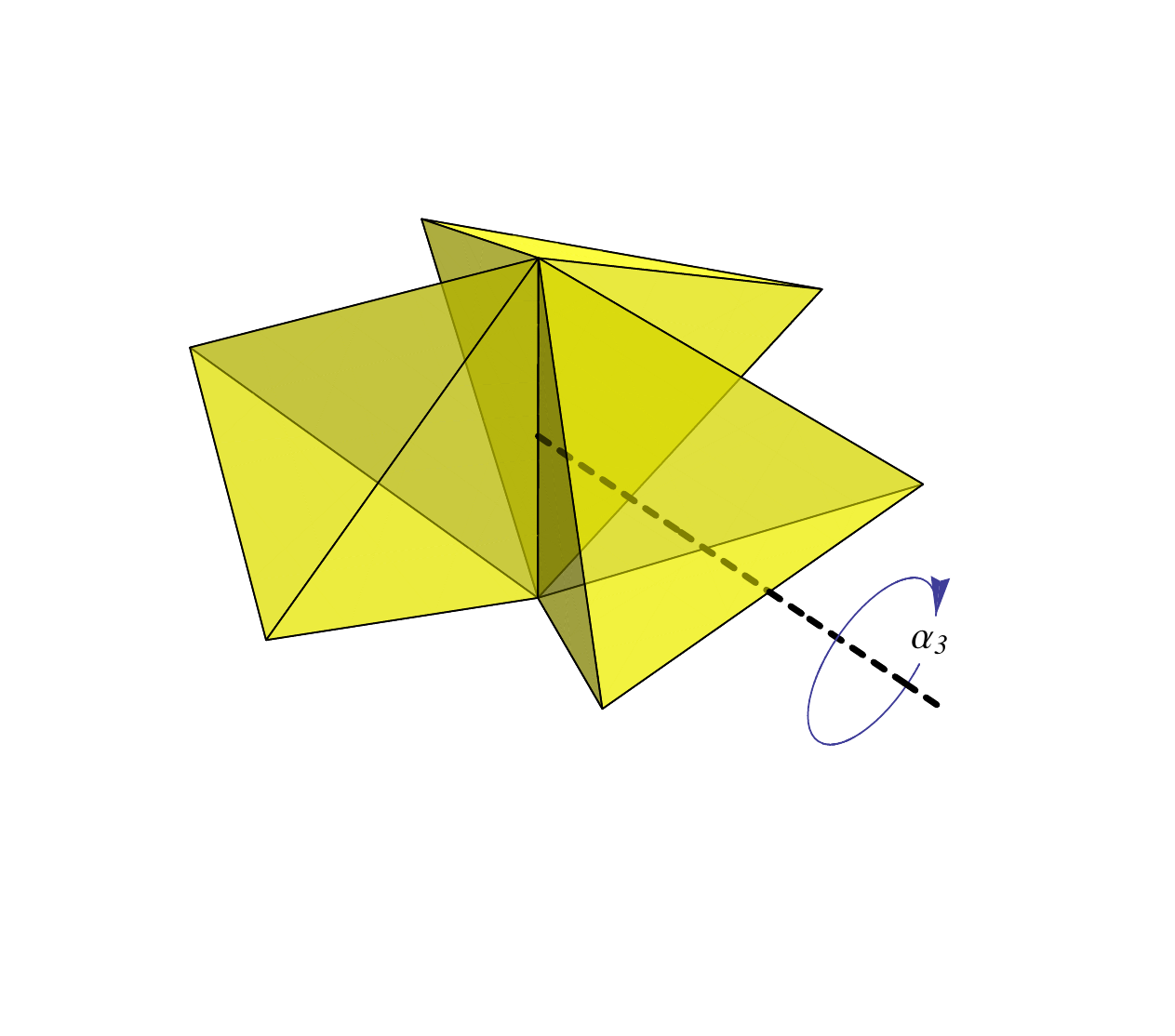}
		\caption{Three tetrahedra arranged about a common edge. In this arrangement, large gaps exist between the faces of adjacent tetrahedra. Each tetrahedron is to be rotated by $\alpha_3$ about an axis passing between the midpoints of its central and peripheral edges.}
		\label{F:tgA}
	\end{subfigure}
	\quad
	\begin{subfigure}[t]{0.3\textwidth}
		\centering
		\includegraphics[width=\textwidth]{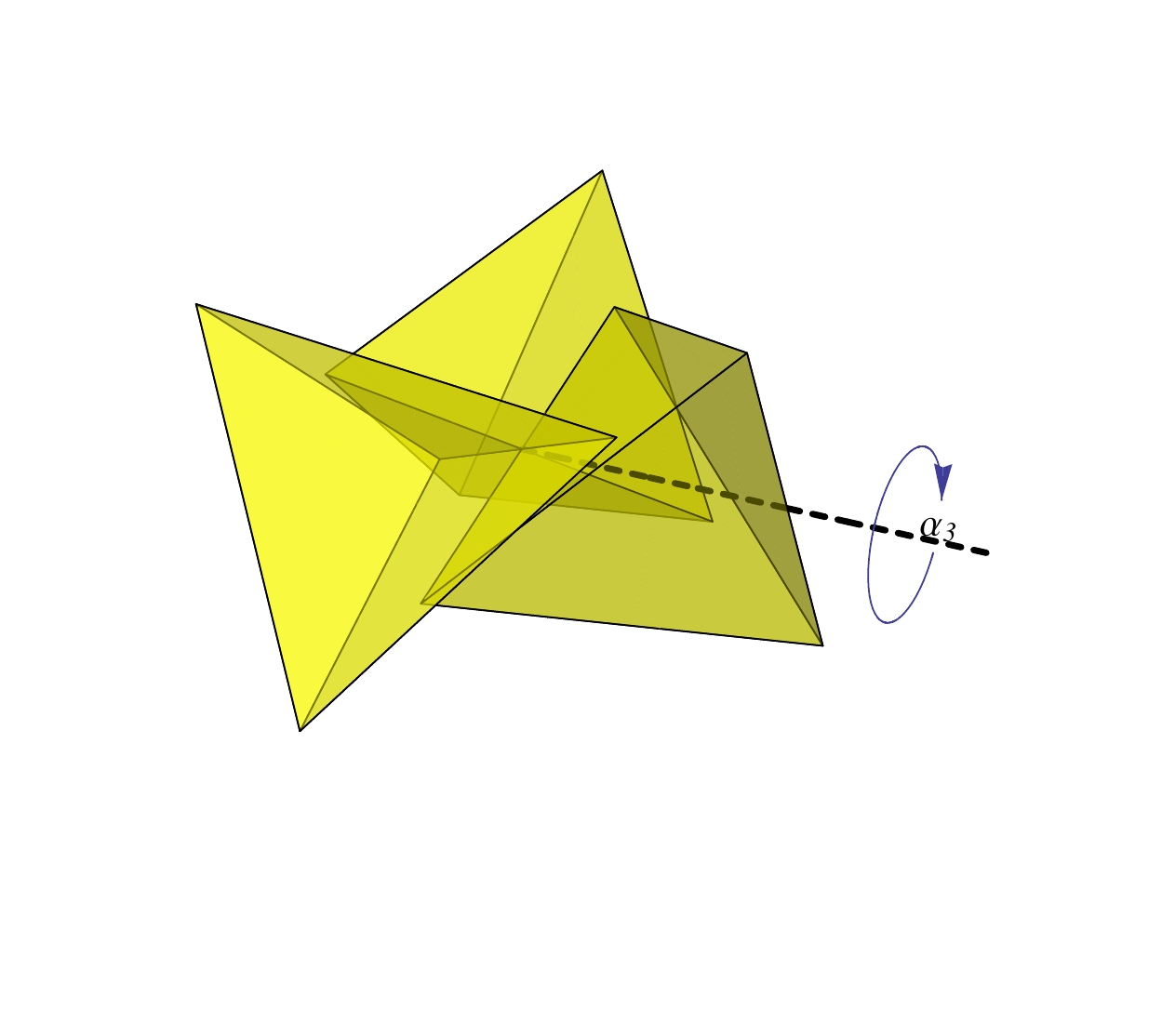}
		\caption{The tetrahedra after rotation. In this arrangement, the faces of adjacent tetrahedra have been brought into contact with one another.}
		\label{F:tgB}
	\end{subfigure}
	\quad
	\begin{subfigure}[t]{0.3\textwidth}
		\centering
		\includegraphics[width=\textwidth]{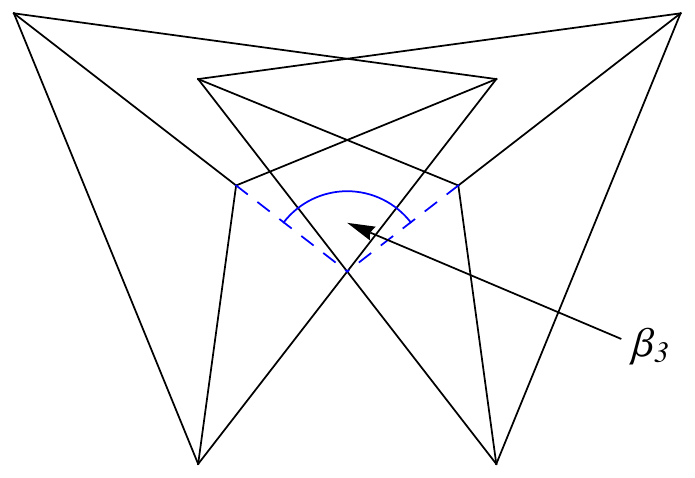}
		\caption{A projection of a ``face junction'' between two coincident faces in \ref{F:tgB}. The angular displacement between the faces is $\beta_3$.}
		\label{F:tgC}
	\end{subfigure}
	\caption{Arranging and rotating three tetrahedra (as done in the $n=5$ case) to ``close up'' gaps between adjacent tetrahedra. When this operation is performed, the angle $\beta_3 = \frac{2\pi}{3} - \beta$ is produced in the projection of a ``face junction.''}
\end{figure}

\begin{figure}[t]
	\label{F:twg}
	\centering
	\begin{subfigure}[t]{0.3\textwidth}
		\centering
		\includegraphics[width=\textwidth]{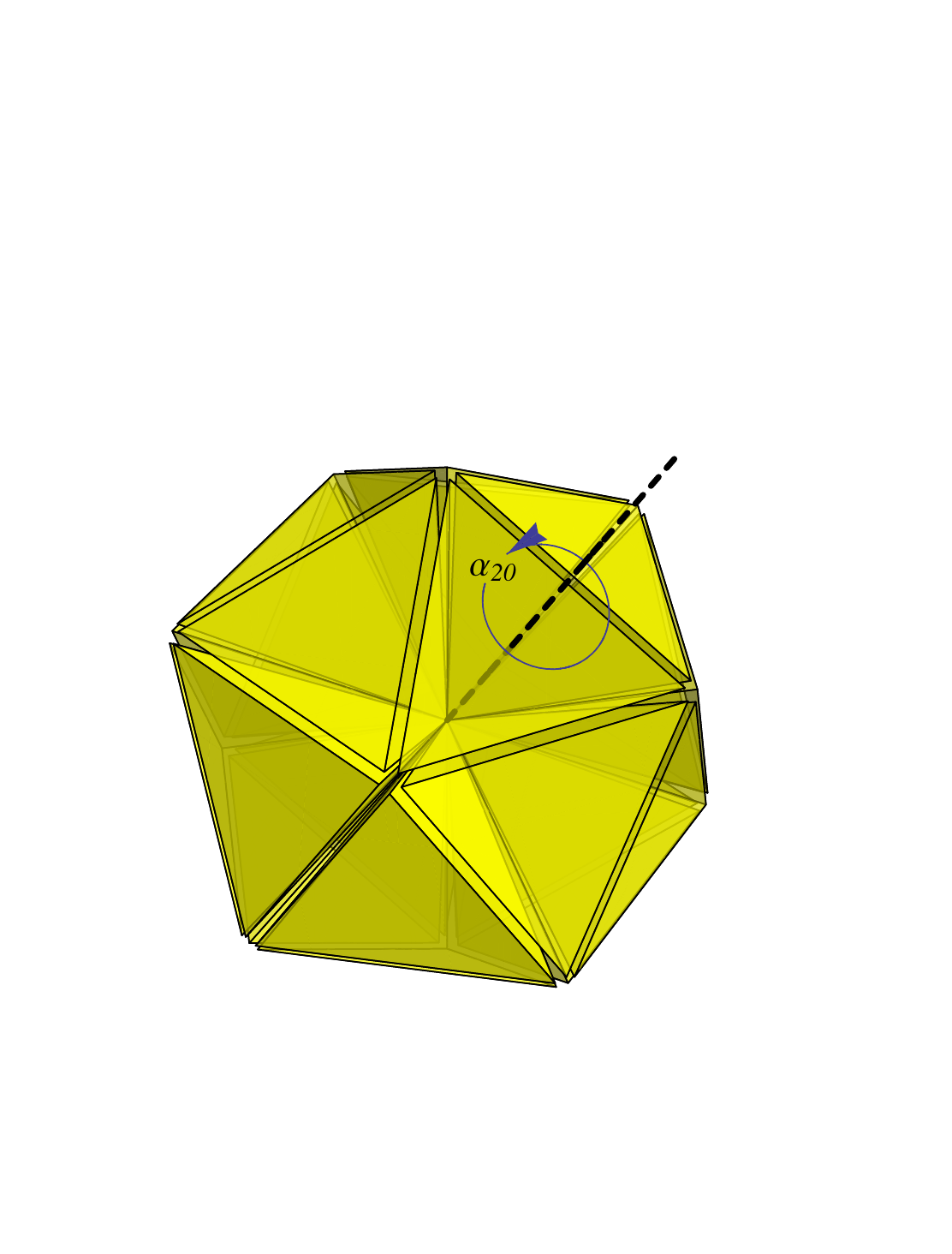}
		\caption{Twenty tetrahedra arranged with icosahedral symmetry about a common central vertex. In this arrangement, gaps exist between faces of adjacent tetrahedra. Each tetrahedron is to be rotated by $\alpha_{20}$ about an axis passing from the central vertex through its exterior face.}
		\label{F:twgA}
	\end{subfigure}
	\quad
	\begin{subfigure}[t]{0.3\textwidth}
		\centering
		\includegraphics[width=\textwidth]{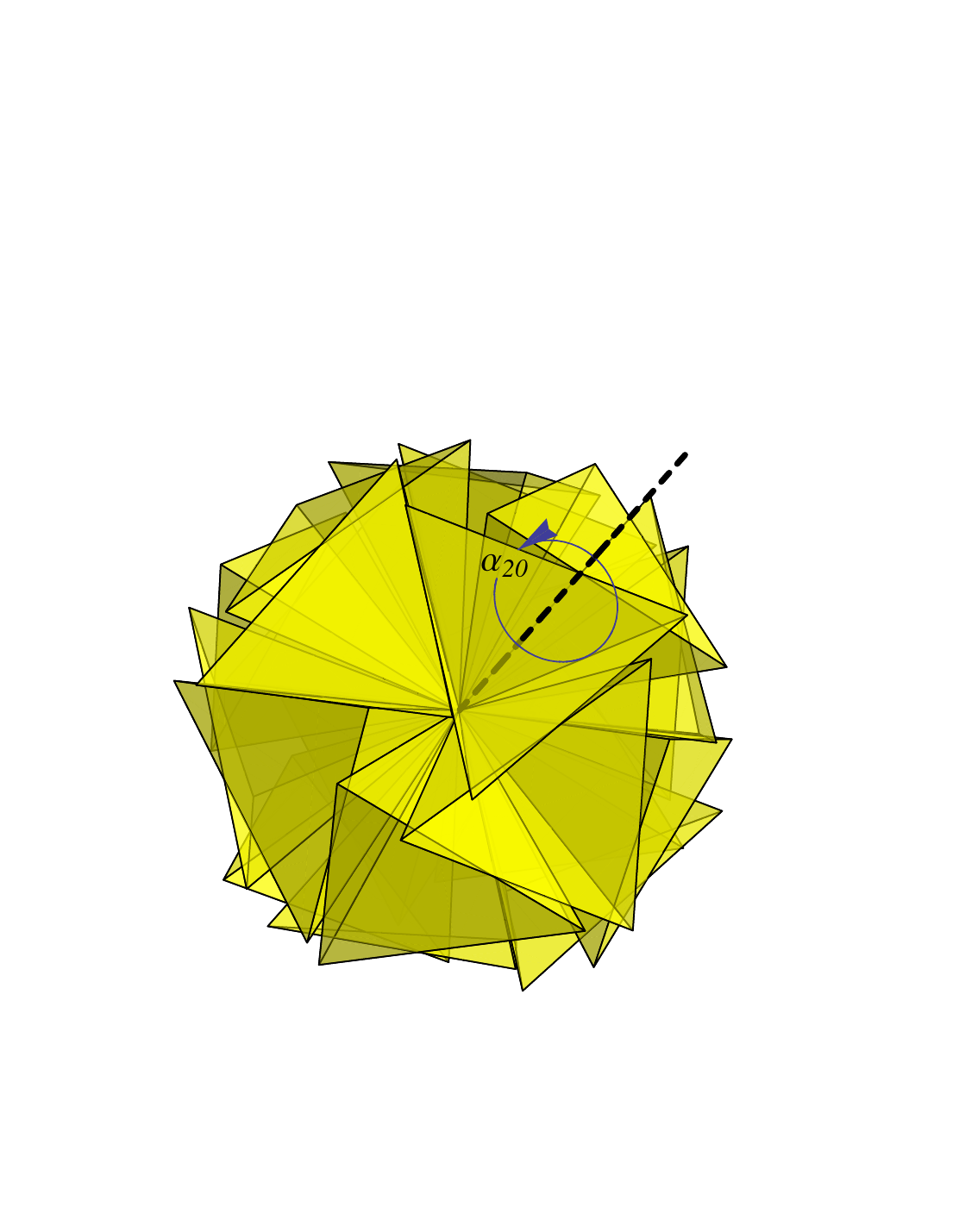}
		\caption{The tetrahedra after rotation. Like in the cases above, face of adjacent tetrahedra have been brought into contact.}
		\label{F:twgB}
	\end{subfigure}
	\quad
	\begin{subfigure}[t]{0.3\textwidth}
		\centering
		\includegraphics[width=\textwidth]{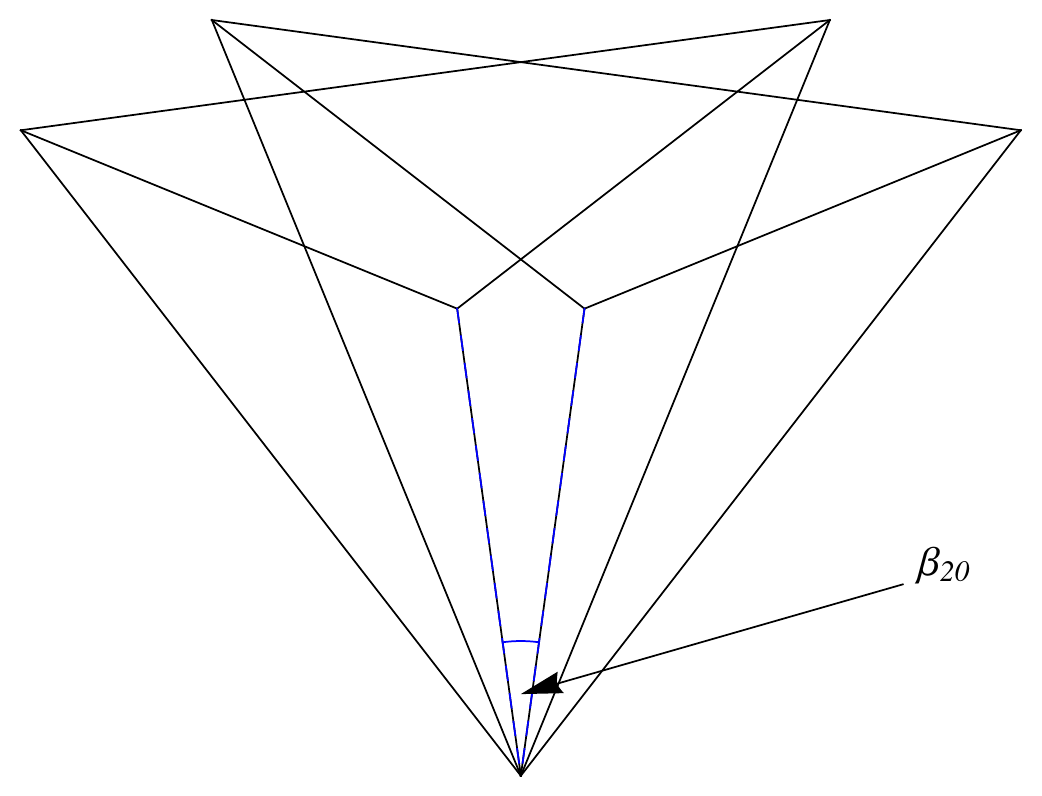}
		\caption{A projection of the ``face junction'' between two coincident faces in \ref{F:twgB}. The angular displacement between the faces is $\beta_{20} = \beta_5 = \beta$.}
		\label{F:twgC}
	\end{subfigure}
	\caption{When 20 tetrahedra are organized into an icosahedral arrangement, gaps between adjacent tetrahedra may be ``closed'' by performing a rotation of each tetrahedron by $\alpha_{20}$ about an axis passing from the central vertex through each tetrahedron's exterior face. When this is done, an angle of $\beta$ is produced in the projection of faces in a ``face junction.''}
\end{figure}

The present document is focused on the angle $\beta = \arccos \left(\left(3\phi-1\right)/4\right)$, which is, in fact, the angle obtained in the ``face junction'' produced by executing the above procedure for $n = 5$ tetrahedra (see Figure~\ref{F:fivegroup}). It is interesting, however, that a simple relationship may be established between this angle and $\beta_3$. By evenly arranging three tetrahedra about an edge and rotating each through the axis extending between the central and peripheral edge midpoints, the face junction depicted in Figure~\ref{F:tgC} is obtained. The angle between faces in this junction, $\beta_3$, may be related to $\beta$ in the following way:
\newline
\begin{equation}
	\frac{2\pi}{3} - \beta = \beta_3.
	\label{beta1}
\end{equation}
\newline
To see this, note that 
\newline
\begin{equation}\label{E:beta2}
	\frac{2\pi}{3} - \arccos \left(\frac{3\phi - 1}{4}\right) = 2 \arctan \left( \frac{\sqrt{\cos^2 \frac{\gamma}{2} - \cos^2 \frac{\theta}{2}}}{\cos \frac{\theta}{2}} \right)
\end{equation}
\newline
gives the solution $\theta = \pm \frac{2\pi}{3}$, which reduces the right hand side of \eqref{E:beta2} to $\beta_3$.

In this section, we have produced two aggregates of tetrahedra whose face junctions bear a relationship to the angle $\beta = \arccos \left(\left( 3 \phi - 1 \right)/4\right)$. In the section that follows, we will locate this angle in the face junctions produced through rotations of tetrahedra about a common vertex.

\begin{figure}[t]
	\centering
	\begin{subfigure}[t]{0.3\textwidth}
		\centering
		\includegraphics[width=\textwidth]{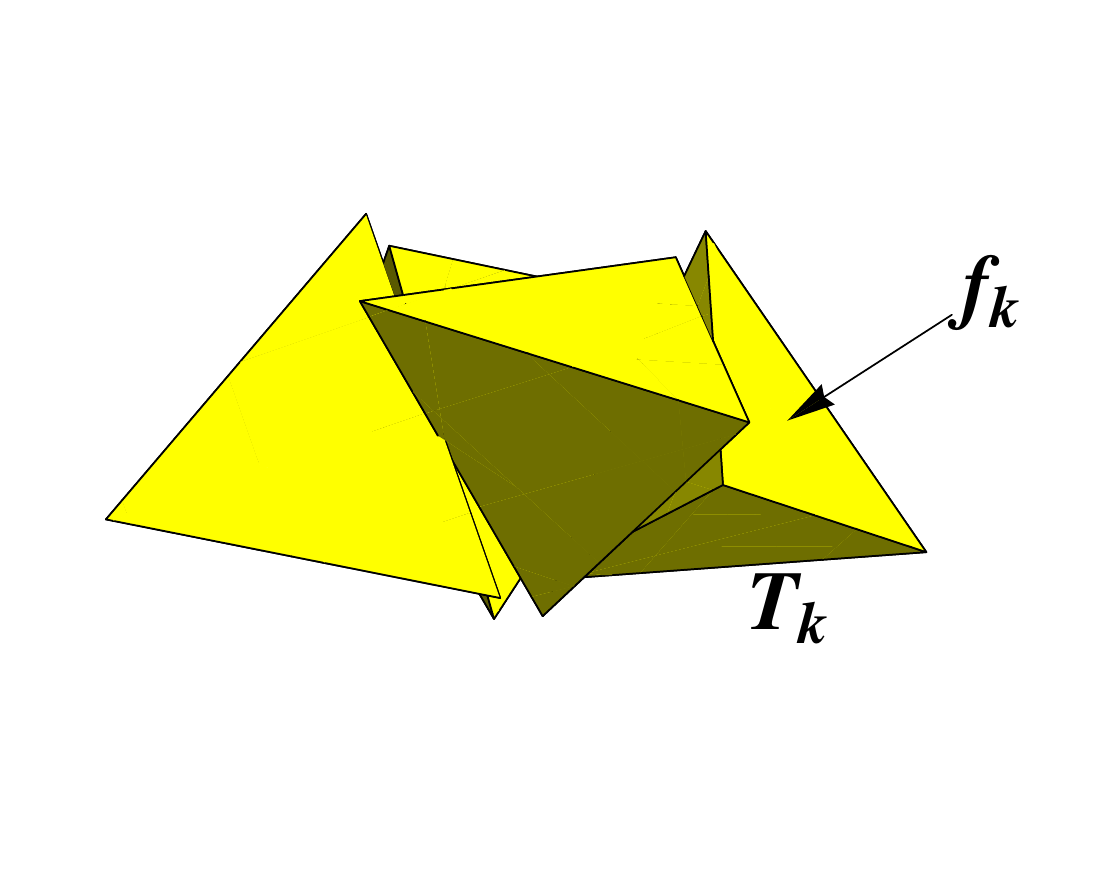}
		\caption{A segment of a modified BC helix with face $f_k$ identified on tetrahedron $T_k$.}
		\label{F:philixA}
	\end{subfigure}
	\quad
	\begin{subfigure}[t]{0.3\textwidth}
		\centering
		\includegraphics[width=\textwidth]{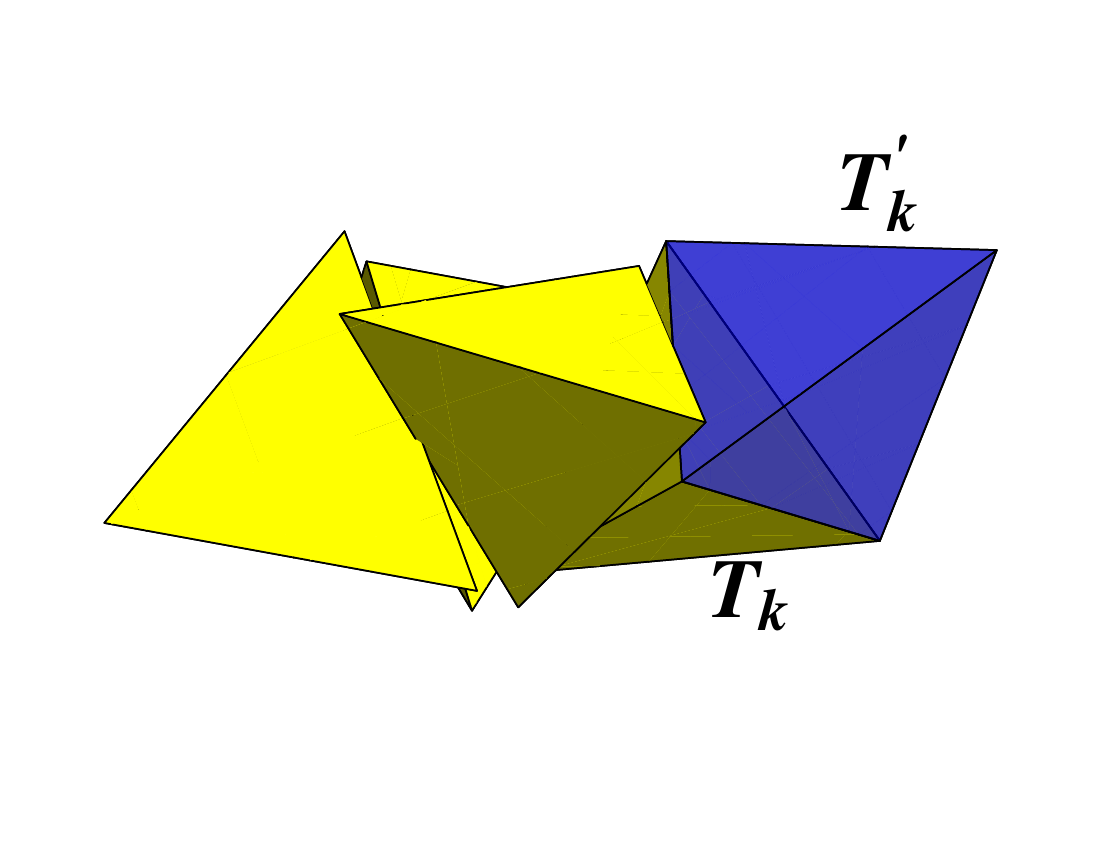}
		\caption{An intermediate tetrahedron, $T_k^\prime$ (shown in blue), is appended (face-to-face) to $f_k$ on $T_k$.}
		\label{F:philixB}
	\end{subfigure}
	\quad
	\begin{subfigure}[t]{0.3\textwidth}
		\centering
		\includegraphics[width=\textwidth]{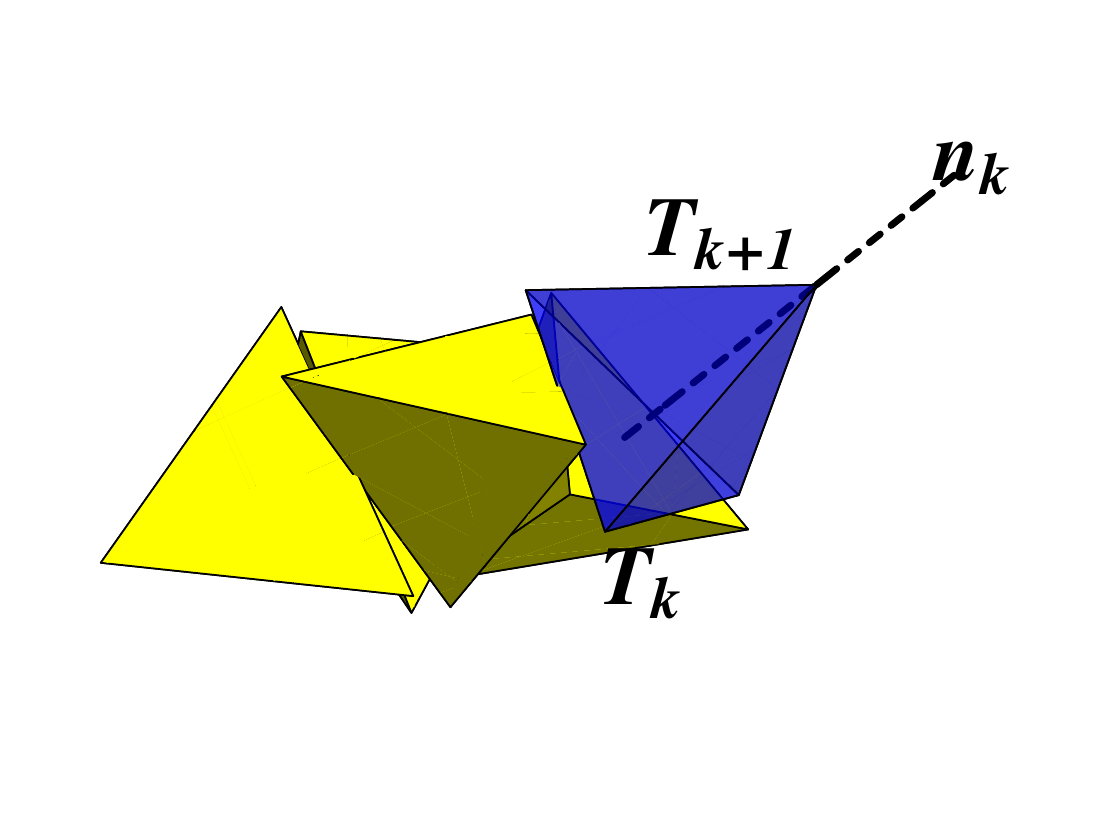}
		\caption{Finally, $T_{k+1}$ is obtained by rotating ${T_k^\prime}$ through the angle $\beta$ about the axis $n_k$.}
		\label{F:philixC}
	\end{subfigure}
	\caption{Assembly of a modified BC helix.}
	\label{F:philix}
\end{figure}

\subsection{Aggregates about a common vertex}\label{S:comvertex}

Consider the icosahedral aggregation of 20 tetrahedra depicted in Figure~\ref{F:twgA}. The face junction of Figure~\ref{F:twgC} is obtained when each tetrahedron is rotated by an angle of
\newline
\begin{equation}
	\alpha_{20} = \arccos \left( \frac{\phi^2}{2\sqrt{2}} \right)
\end{equation}
\newline
about an axis extending between the center of its exterior face and the arrangement's central vertex. As above, this operation ``closes'' gaps between tetrahedra by bringing adjacent faces into contact. Interestingly, the face junction obtained here consists of tetrahedra with a rotational displacement equal to the one obtained in the case of five tetrahedra arranged about a common edge above, i.e., $\beta_{20} = \beta_5 = \beta$. (It should be noted, however, that $\alpha_{20}$ and $\beta_{20}$ are not produced by Equations~\eqref{E:alphan} and \eqref{E:betan}, respectively, as those formul\ae\ are only valid for $3 \leq n \leq 5$.)

In all of the cases described above, gaps are ``closed'' and ``junctions'' are produced between adjacent tetrahedra in such a way that the angle $\beta = \arccos \left(\left(3\phi-1\right)/4\right)$ appears in some fashion in the angular displacement between coincident faces. For the case of 5 tetrahedra about a central edge and 20 tetrahedra about a common vertex, this angle is observed directly. For the case of 3 tetrahedra about a central edge, the angular displacement between faces is closely related: $\beta_3 = \frac{2\pi}{3} - \beta$. 

We now turn to an arrangement obtained by directly imparting an angular displacement of $\beta$ between adjacent pairs of tetrahedra in a linear, helical fashion known as the Boerdijk-Coxeter helix. An interesting result of performing this action is that a previously aperiodic structure is transformed into one with translational and rotational symmetries.

\begin{figure}[t]
	\centering
	\begin{subfigure}[t]{0.3\textwidth}
		\centering
		\includegraphics[width=\textwidth]{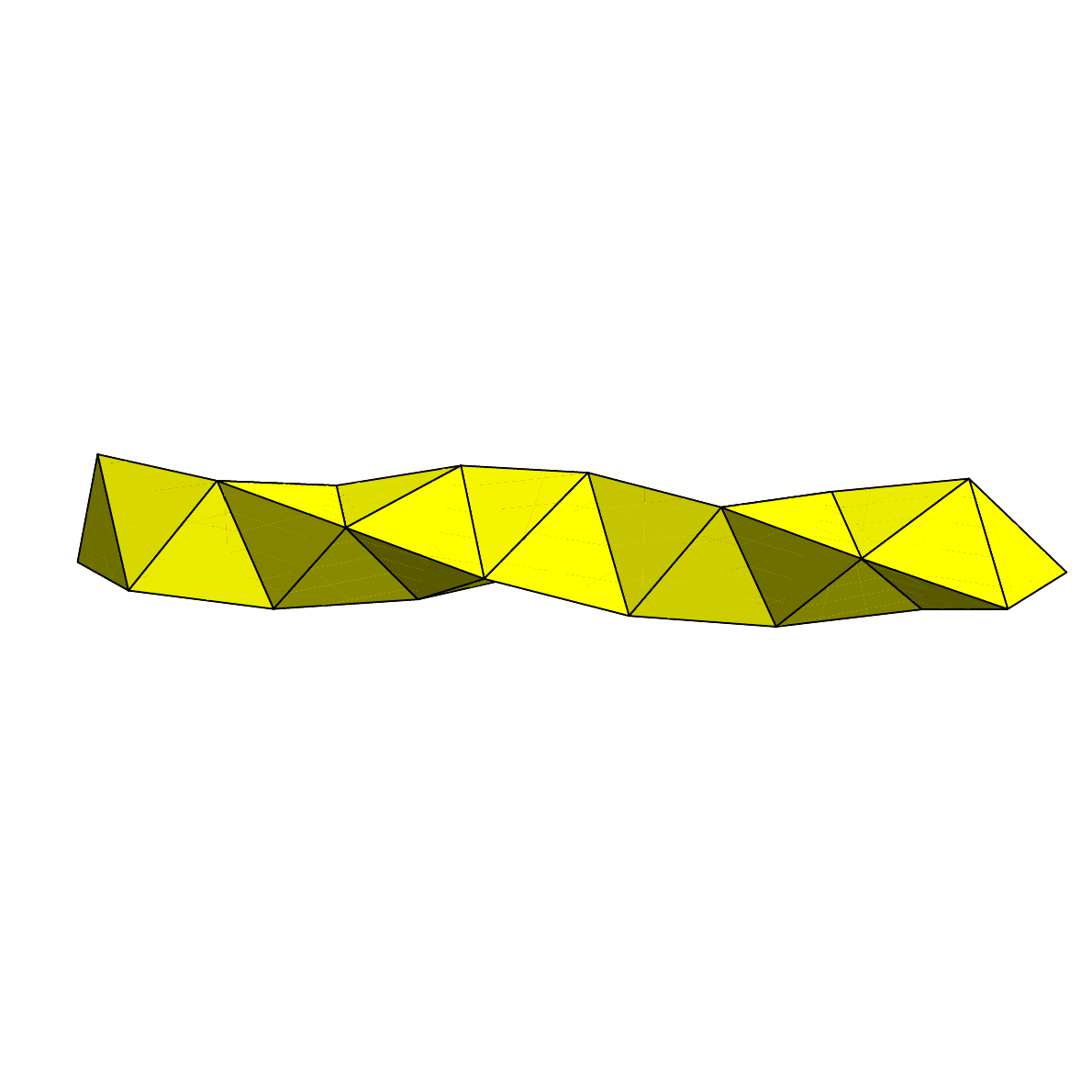}
		\caption{A right-handed BC helix.}
		\label{F:bchelix}
	\end{subfigure}
	\quad
	\begin{subfigure}[t]{0.3\textwidth}
		\centering
		\includegraphics[width=\textwidth]{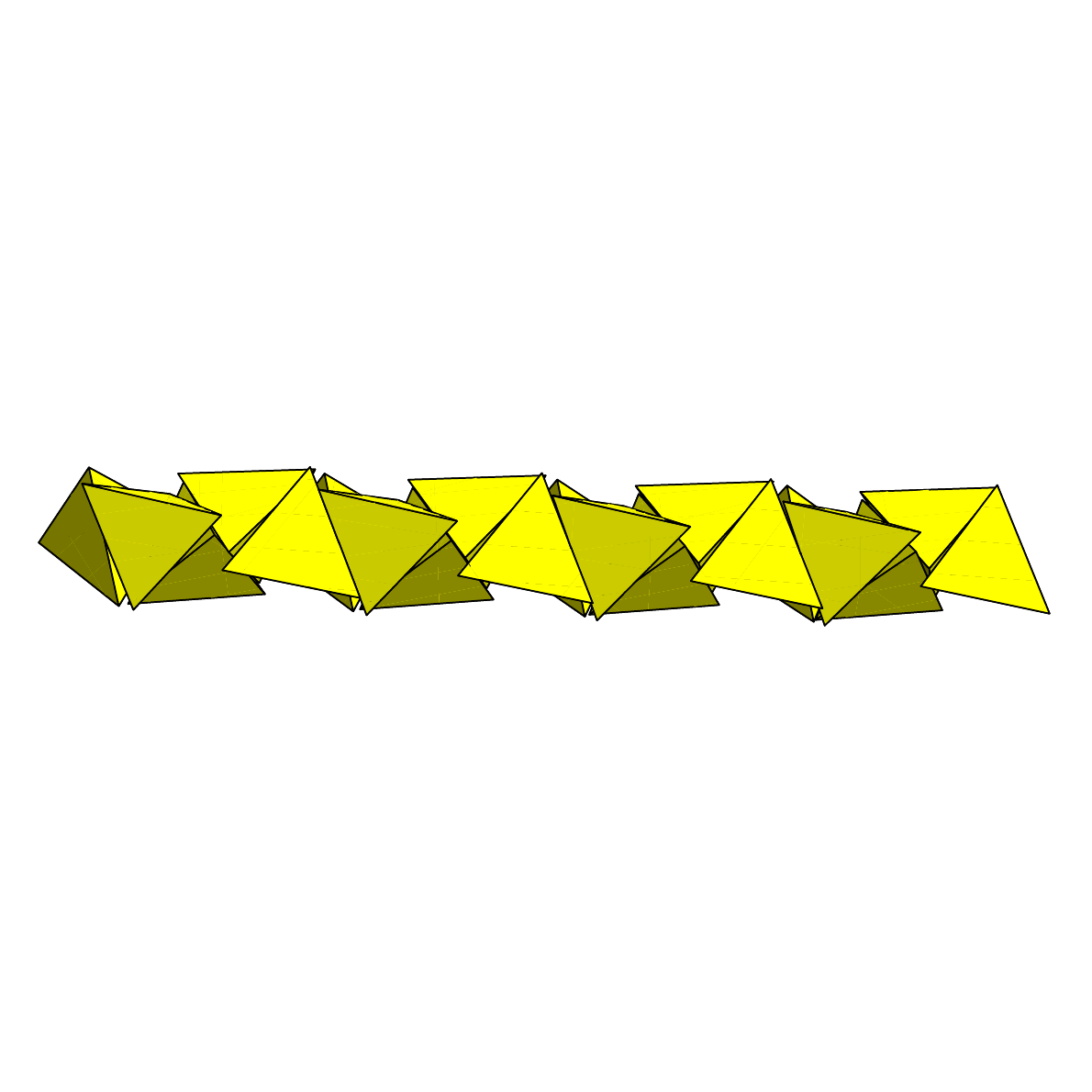}
		\caption{A ``5-BC helix'' may be obtained by appending and rotating tetrahedra by $\beta$ using the same chirality of the underlying helix.}
		\label{F:5bchelix}
	\end{subfigure}
	\quad
	\begin{subfigure}[t]{0.3\textwidth}
		\centering
		\includegraphics[width=\textwidth]{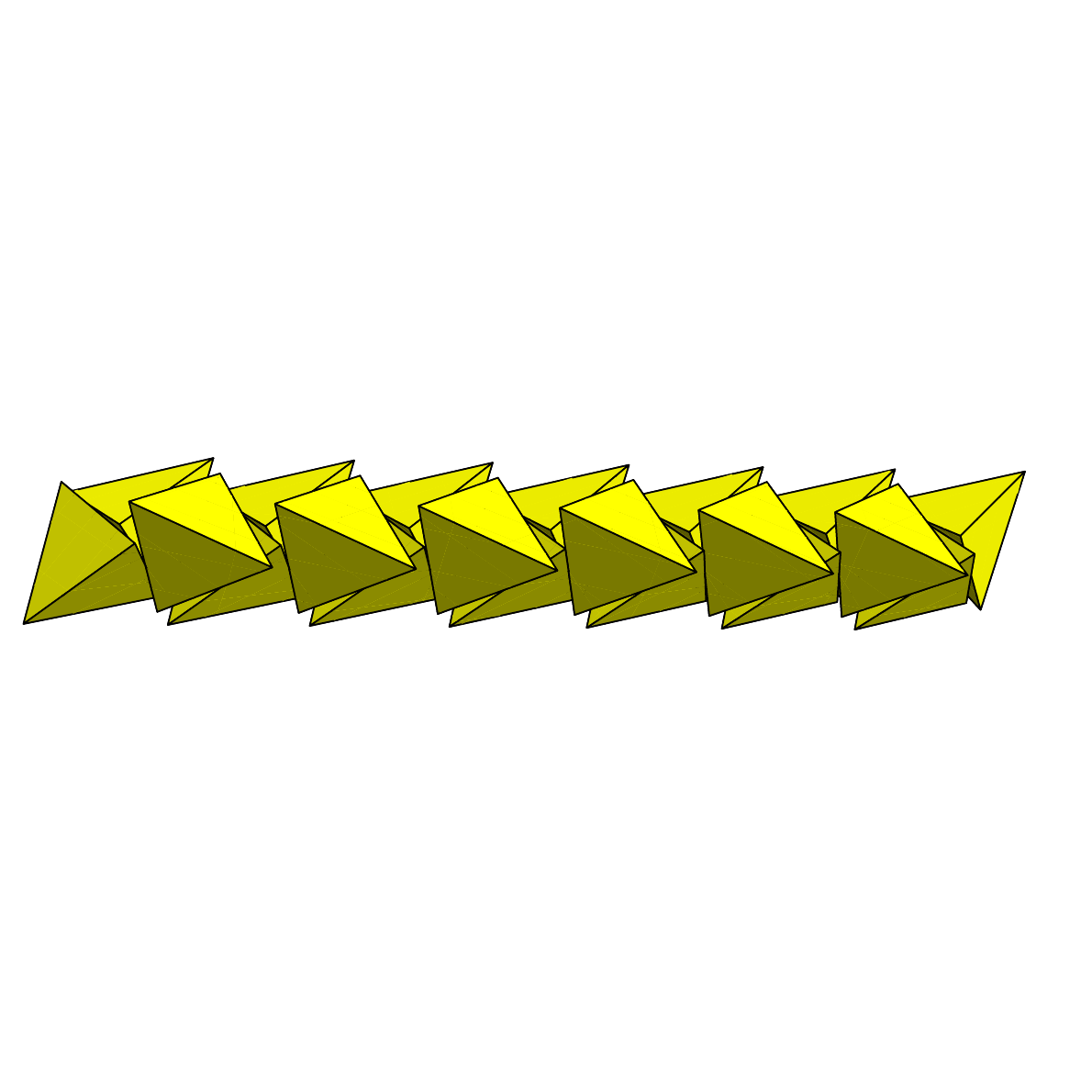}
		\caption{A ``3-BC helix'' may be obtained by appending and rotating tetrahedra by $\beta$ using the opposite chirality of the underlying helix.}
		\label{F:3bchelix}
	\end{subfigure}
	\caption{Canonical and modified Boerdijk-Coxeter helices.}
	\label{F:mbchelices}
\end{figure}

\subsection{Periodic, helical aggregates}\label{S:helix}

In Sections~\ref{S:comedge} and \ref{S:comvertex} we described a procedure by which initial arrangements of tetrahedra were transformed so that adjacent pairs of tetrahedra were brought together to touch. In each of these structures, coincident faces are displaced by an angle equal or closely related to $\beta$. Here, will will construct two periodic, helical chains of tetrahedra by directly inserting an angular offset by $\beta$ between each successive member of the chain. For their close relationship with the Boerdijk-Coxeter helix, we refer to these structures by the term \emph{modified BC helices}.

\begin{figure}[t]
	\centering
	\begin{subfigure}[t]{0.3\textwidth}
		\centering
		\includegraphics[width=\textwidth]{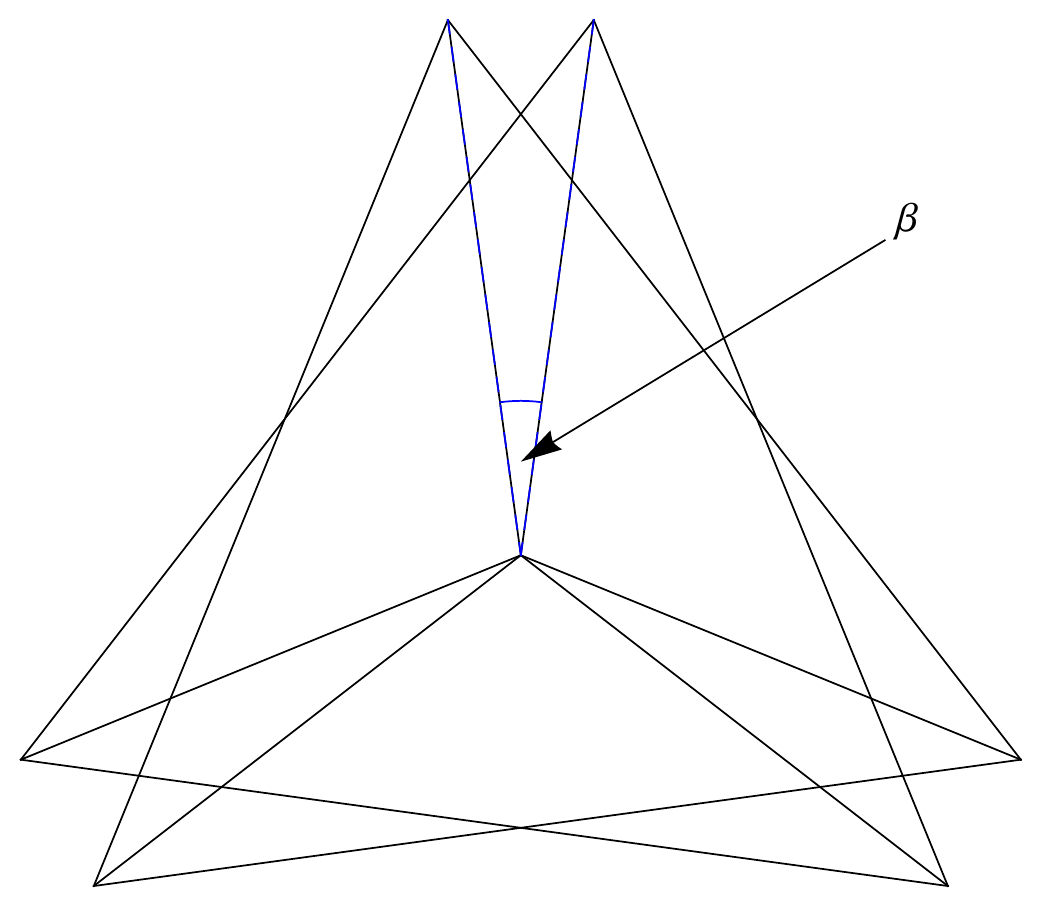}
		\caption{A ``face junction'' between tetrahedra of a 3-- or 5--BC helix.}
		\label{F:bcjunction}
	\end{subfigure}
	\quad
	\begin{subfigure}[t]{0.3\textwidth}
		\centering
		\includegraphics[width=\textwidth]{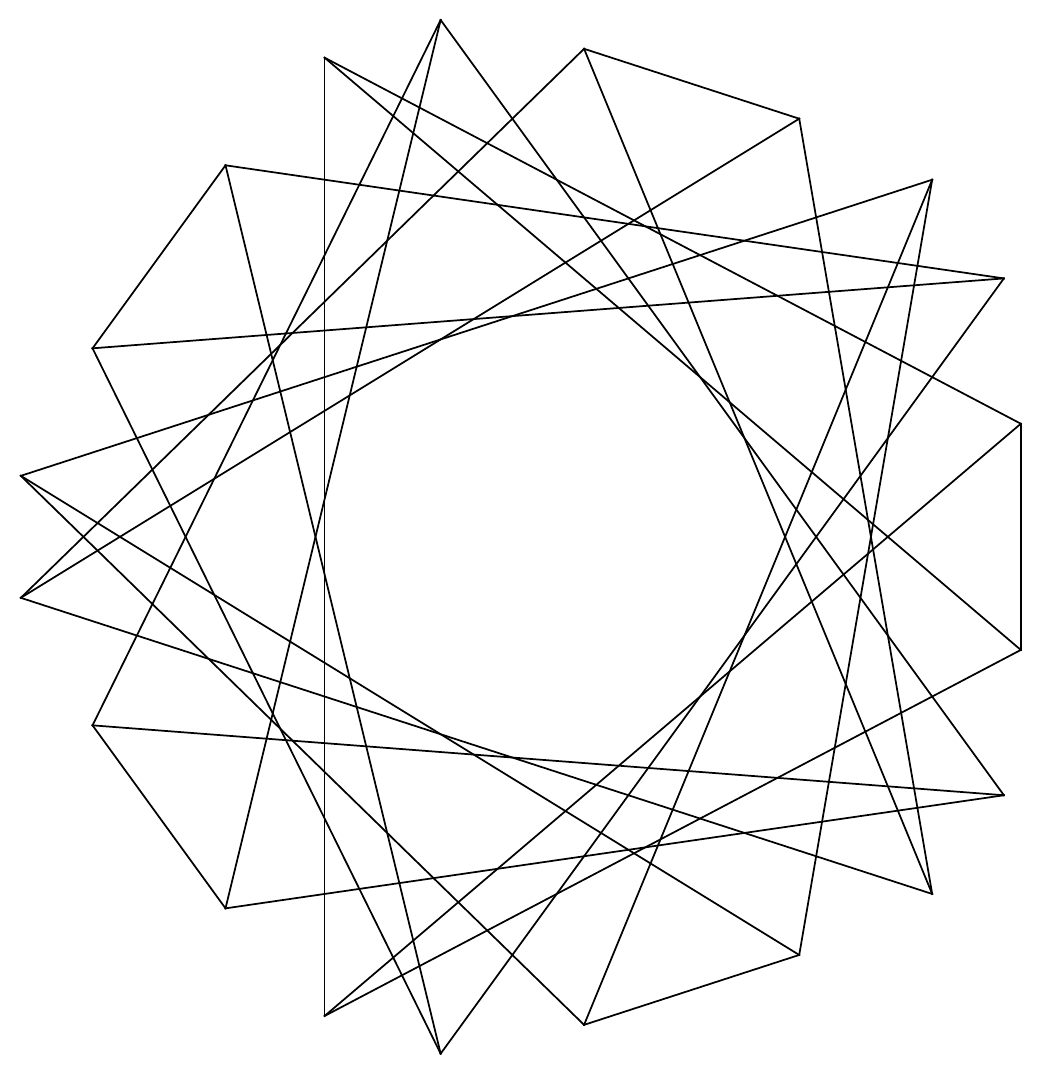}
		\caption{A projection of the 5--BC helix along its central axis showing five-fold symmetry.}
		\label{F:5bcproj}
	\end{subfigure}
	\quad
	\begin{subfigure}[t]{0.3\textwidth}
		\centering
		\includegraphics[width=\textwidth]{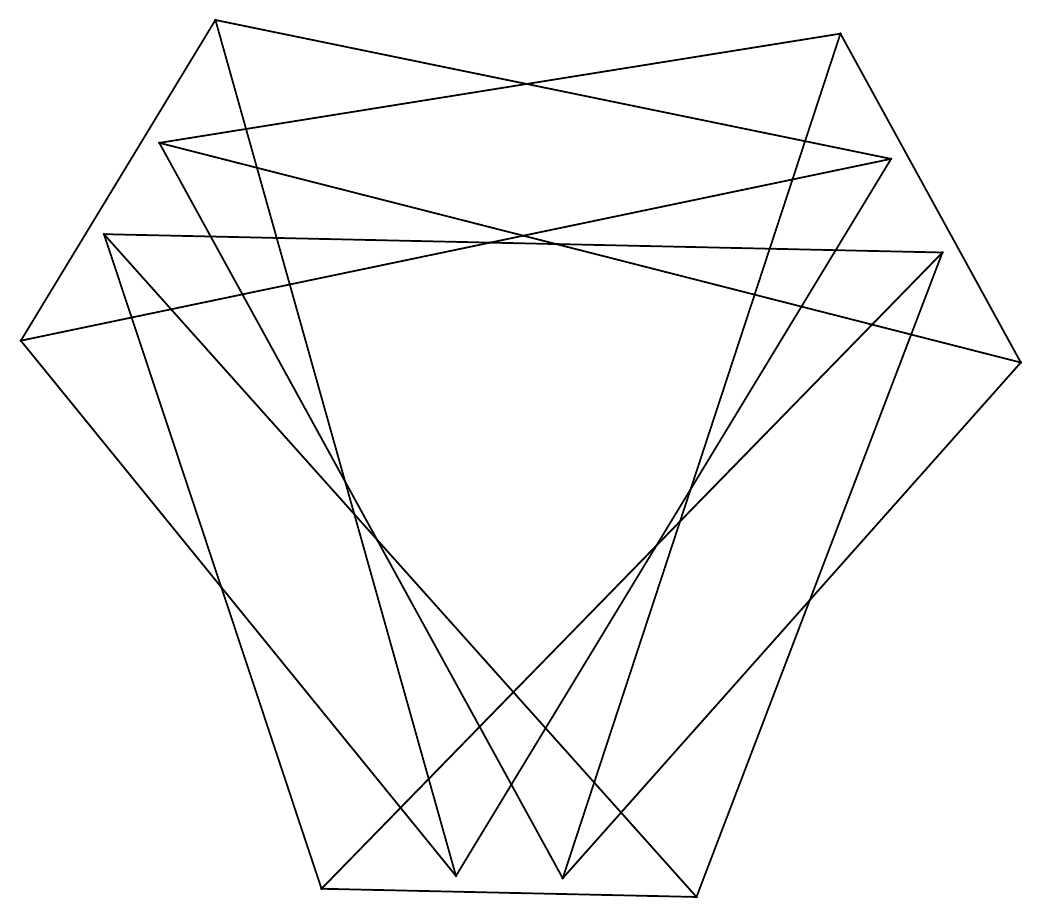}
		\caption{A projection of the 3--BC helix along its central axis showing three-fold symmetry.}
		\label{F:3bcproj}
	\end{subfigure}
	\caption{BC helix projections and face junction.}
	\label{F:bcjproj}
\end{figure}

The construction of a modified BC helix is depicted in Figure~\ref{F:philix}. Starting from a tetrahedron $T_k = \left( v_{k0}, v_{k1}, v_{k2}, v_{k3} \right)$, a face $f_k$ is selected onto which an interim tetrahedron, $T_k^\prime$, is appended. The $\left(k+1\right)^\text{th}$ tetrahedron is obtained by rotating $T_k^\prime$ through an angle $\beta$ about an axis $n_k$ normal to $f_k$, passing through the centroid of $T_k^\prime$. (Note that this automatically produces an angular displacement of $\beta$ between two faces in a ``junction,'' see Figure~\ref{F:bcjunction}.)

The structure that results from this process depends on the sequence of faces $F=\left(f_0,f_1,\ldots,f_k\right)$ selected in order to construct the helical chain. This sequence determines an \emph{underlying chirality} of the helix---i.e., the chirality of the helix formed by the tetrahedral centroids---and plays a pivotal role in the determination of the structure's eventual symmetry. (However, it should be noted, of course, that some sequences of faces do not result in helical structures. Faces cannot be chosen arbitrarily or randomly; they must be selected so as to build a helix.)

By performing the procedure depicted in Figure~\ref{F:philix}, using an angular displacement of $\beta$ between successive tetrahedra, periodic structures are obtained with 3-- or 5--fold symmetry (upon their projections, see Figures~\ref{F:5bchelix} and \ref{F:3bchelix}), depending on the relative chiralities between the rotational displacement and the underlying helix: when \textbf{like} chiralities are used one obtains 5--fold symmetry; when \textbf{unlike} chiralities are used one obtains 3--fold symmetry. In addition to rotational symmetry, these structures are given a linear period, which we quantify here as the number of appended tetrahedra necessary to return to an initial angular position on the helix. For a modified BC helix with a period of \emph{m} tetrahedra, we use the term \emph{m}-BC helix. Accordingly, the procedure described above produces 3-- and 5--BC helices, which are shown in Figure~\ref{F:mbchelices}. (See \cite{sadler2013} for a proof of these structures' symmetries and periodicities.)

\begin{table}[t]
	\centering
	\begin{tabular}{ c c c c c }
		\toprule
		 & & &\multicolumn{2}{c}{Plane Classes}\\
		 \cmidrule(l){4-5}
		Structure & $\alpha_n$ & $\beta_n$ & Before & After\\
		\midrule
		3 tetrahedra, common edge & $\arccos \left( \frac{1}{\sqrt{6}} \right)$ & $\frac{2\pi}{3} - \arccos \left(\frac{3\phi-1}{4}\right)$ & 12 & 9 \\
		4 tetrahedra, common edge & $\frac{\pi}{4}$ & $\frac{\pi}{3}$ & 16 & 4\\
		5 tetrahedra, common edge & $\arccos \left( \frac{\phi^2}{\sqrt{2\left(\phi + 2 \right)}}\right)$ & $\arccos \left(\frac{3\phi-1}{4}\right)$ & 20 & 10\\
		20 tetrahedra, common vertex & $\arccos \left(\frac{\phi^2}{2\sqrt{2}}\right)$ & $\arccos \left(\frac{3\phi-1}{4}\right)$ & 60 & 10\\
		$n$-tetrahedron 3-BC helix & n/a & $\arccos \left(\frac{3\phi-1}{4}\right)$ & $3n+1$ & 9\\
		$n$-tetrahedron 5-BC helix & n/a & $\arccos \left(\frac{3\phi-1}{4}\right)$ & $3n+1$ & 10\\
		\bottomrule
	\end{tabular}
	\caption{Plane class numbers for aggregates described in Section~\ref{S:aggregates}.}
	\label{T:planeclasses}
\end{table}

\section{Curiosities}

We have seen the construction of several aggregates of tetrahedra. Each of these structures contains tetrahedra with coincident faces, offset angularly by $\beta$ or a closely related angle ($\frac{2\pi}{3}-\beta$). We will now explore some of the interesting features of these structures. As already noted, it is interesting that $\beta$ appears in the face junctions of structures generated by ``closing'' gaps between tetrahedral aggregates. It is additionally interesting that when this angle is employed in the construction of a helical chain of tetrahedra, a periodic structure emerges (whereas the canonical BC helix has no non-trivial translational or rotational symmetries). 

The features we will highlight in this section involve the reduction in the overall number of ``plane classes'' and the linear displacements between facial centers in a face function. Here, we say that two planes belong to the same \emph{plane class} if and only if their normal vectors are parallel. The \emph{number of plane classes} for a collection of tetrahedra, then, is defined as the number of distinct plane classes comprising the collection's two-dimensional faces. By rotating tetrahedra so as to bring faces into contact (as in Sections~\ref{S:comedge} and \ref{S:comvertex}), or rotating tetrahedra to obtain periodicity (as in Section~\ref{S:helix}) the overall number of plane classes for an aggregate is reduced. Clearly, as the values of $\alpha_3$, $\alpha_5$, and $\alpha_{20}$ are such that they bring faces of adjacent tetrahedra into contact, we would expect to see a reduction in the number of plane classes in the corresponding aggregations of tetrahedra featured in Sections~\ref{S:comedge} and \ref{S:comvertex}. It is interesting, however, that rotation of the tetrahedra in a BC helix by $\beta$ (observed in the face junctions of Figures~\ref{F:fgC}, \ref{F:tgC}, and \ref{F:twgC}) obtains a reduction from an arbitrarily large number of plane classes ($3 n + 1$, where $n$ is the number of tetrahedra in the helix) to relatively small numbers: 9 plane classes in the case of the 3-BC helix, 10 plane classes in the case of the 5-BC helix. Table~\ref{T:planeclasses} provides the numbers of plane classes for the tetrahedral aggregates described in Section~\ref{S:aggregates} before and after their transformations.

\begin{figure}[t]
	\centering
	\includegraphics[width=\textwidth]{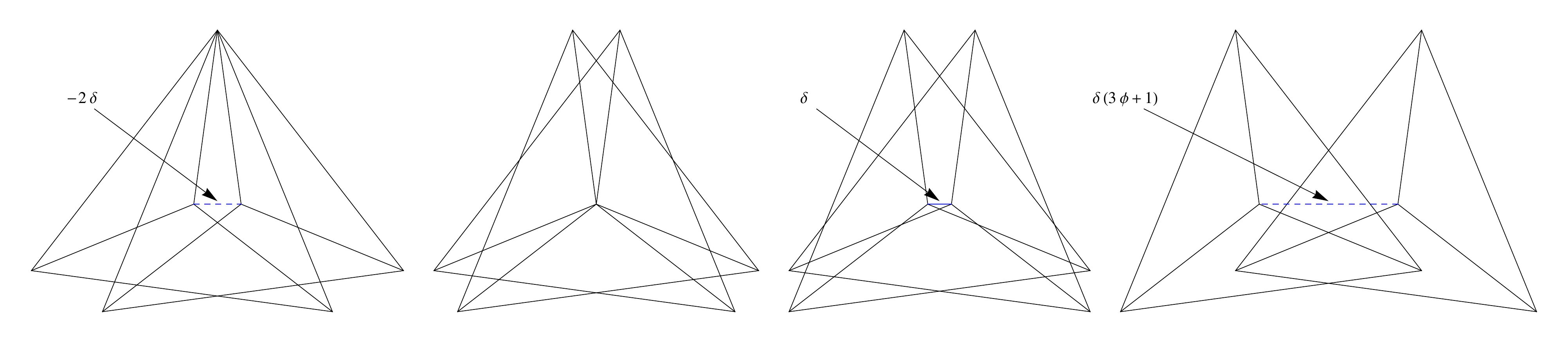}
	\caption{Side-by-side comparison of face junctions. All face junctions may be obtained by translation of the (projected) tetrahedra of the 3-- and 5-BC helix (pictured center left) by integer and golden ratio-based multiples of $\delta = \left<\frac{a}{2\phi^2\sqrt{6}},0\right>$, where $a$ is the tetrahedron edge length, the displacement between tetrahedra of a face junction for 5 tetrahedra about a common edge (pictured center right). From left to right the displacements between the tetrahedra are $-2\delta$, $0$, $\delta$, and $\left(3\phi+1\right)\delta$.}
	\label{F:facejunctions}
\end{figure}

Finally, an appealing feature is observed in the face junction projections of the tetrahedral aggregates discussed in this paper. Figure~\ref{F:facejunctions} provides a side-by-side comparison of these face junctions. As the angular displacement between tetrahedra in all face junctions is related to $\beta$, we can see that translation of a tetrahedron in one junction can produce any of the other junctions. Let the displacement between tetrahedra in the face junction of 5 tetrahedra about a common edge be denoted by $\delta = \left<\frac{a}{2\phi^2 \sqrt{6}}, 0\right>$, where $a$ is the tetrahedron edge length. Starting from the face junction of a 3-- or 5-BC helix, the remaining face junctions corresponding to 20 tetrahedra about a vertex, 5 tetrahedra about an edge, and 3 tetrahedra about an edge may be obtained by translating a tetrahedron of the junction by $-2\delta$, $\delta$, and $\left(3\phi+1\right)\delta$, respectively.

\section{Conclusion}

In this paper we have presented the construction of several aggregates of tetrahedra. In each case, the construction process involved rotations of tetrahedra by a value related to $\beta=\arccos\left(3\phi-1\right)/4$. The structures produced here have several notable features: faces of tetrahedra are made to touch (``closing'' previously existing gaps between tetrahedra), aperiodic structures are imparted with periodicity, and the total number of plane classes is reduced. The purpose of the present document, however, is not merely descriptive; it is hoped that these notable features have generated interest in the reader of rotational transformations of tetrahedra involving the angle $\beta$. In particular, it is desired that further observations may be found which transform aggregates of tetrahedra such that faces are brought into contact and the number of plane classes is reduced.

\end{document}